\documentclass[12pt]{amsart}

\usepackage{amsmath,amssymb,amscd,verbatim}

\newcommand{\mc}{\mathcal}

\newcommand{\sub}{\subseteq}

\newcommand{\ol}{\overline}

\newcommand{\lra}{\Leftrightarrow}

\newcommand{\ra}{\Rightarrow}

\newcommand{\sm}{\setminus}

\newcommand{\al}{\alpha}

\newcommand{\be}{\beta}

\newcommand{\La}{\Lambda}

\newcommand{\Ga}{\Gamma}

\DeclareMathOperator{\Spec}{-Spec}

\DeclareMathOperator{\Max}{-Max}

\newcommand{\tmax}{t\text{-Max}}

\newcommand{\rx}{R[\{X_{\al}\}]}

\newcommand{\op}{\operatorname}

\newcommand{\rxx}{R(\{X_{\al}\})}

\newcommand{\rxxx}{R\langle \{X_{\al}\} \rangle}

\newtheorem{theorem}{Theorem}[section]

\newtheorem{lemma}[theorem]{Lemma}

\newtheorem{proposition}[theorem]{Proposition}

\newtheorem{corollary}[theorem]{Corollary}

\theoremstyle{definition}

\newtheorem{definition}[theorem]{Definition}

\newtheorem{remark}[theorem]{Remark}

\newtheorem{example}[theorem]{Example}

\begin{document}

\title{The $t\#$-property for integral domains}

\author{Stefania Gabelli}

\address{Dipartimento di Matematica\\ Universit\`{a} degli Studi Roma Tre\\
 Largo San L. Murialdo, 1\\
00146 Roma, Italy}

\email{gabelli@mat.uniroma3.it}

\author{Evan Houston}

\address[Evan Houston and Thomas G. Lucas]{Department of Mathematics\\
 University of North Carolina at Charlotte\\
 Charlotte, NC 28223 U.S.A.}

\email[Evan Houston]{eghousto@email.uncc.edu}

\author{Thomas G. Lucas}
\email[Thomas G. Lucas]{tglucas@email.uncc.edu}

\thanks{All three authors acknowledge support from the Cultural
Co-operation Agreement between Universit\`{a} degli Studi Roma Tre and
The University of North Carolina at Charlotte.}

\maketitle

\date{}




\section*{introduction}

Let $R$ denote an integral domain with quotient field $K$.  Then
$R$ is said to be a \emph{$\#$-domain} or to satisfy the
\emph{$\#$-condition} if $\bigcap_{M \in \mc M_1}R_M \ne
\bigcap_{M \in \mc M_2}R_M$ whenever $\mc M_1$ and $\mc M_2$ are
distinct subsets of the set of maximal ideals of $R$.  Pr\"ufer
domains satisfying the $\#$-condition were first studied in
\cite{g3} and \cite{gh}.  Domains each of whose overrings satisfy
the $\#$-condition were also studied in \cite{gh} (in the Pr\"ufer
case); these domains have come to be called \emph{$\#\#$-domains}.

Although the papers mentioned above contain very interesting
results, those results are essentially restricted to the class of
Pr\"ufer domains.  This paper represents an effort to extend, by a
modification of the definitions, results about the $\#$- and
$\#\#$-conditions to a much wider class of domains.  In the first
section, we introduce the $t\#$-condition: A domain $R$ satisfies
the $t\#$-condition if $\bigcap_{M \in \mc M_1} R_M \ne \bigcap_{M
\in \mc M_2} R_M$ for any two distinct subsets $\mc M_1, \mc M_2$
of the set of maximal $t$-ideals of $R$. We discuss the extent to
which the properties shown in \cite{gh} to be equivalent to the
$\#$-property carry over to our setting. For example,
\cite[Theorem 1 (a) $\lra$ (b)]{gh} states that the domain $R$ has
the $\#$-property if and only if each maximal ideal $M$ of $R$
contains a finitely generated ideal which is contained in no other
maximal ideal of $R$; we show that this result has a natural
counterpart in the class of $v$-coherent domains (which includes
all Noetherian domains). (All relevant definitions are given
below.)  In addition, we show that for \emph{any} domain $R$, $R$
has the $t\#$-property if and only if each maximal $t$-ideal $M$
of $R$ contains a \emph{divisorial} ideal contained in no other
maximal $t$-ideal of $R$. We also give examples to show that
``divisorial'' cannot be replaced by ``finitely generated'' in
general.

In Section 2, we attempt to generalize the $\#\#$-property. In the
case of Pr\"ufer domains, the definition of the $\#\#$-property is
reasonable since the overrings have nice properties (e.g., they
are flat).  To obtain a useful definition of the $t\#\#$-property
for more general classes of rings, however, one must decide which
overrings should be required to have the $t\#$-property. For
example, we could say that $R$ has the $t\#\#$-property if each
$t$-linked overring of $R$ has the $t\#$-property.  Another
possibility is to require that the overrings of $R$ which are
generalized rings of quotients of $R$ should have the
$t\#$-property.  In the end we avoid making a definition at all.
Instead, we explore several classes of overrings, primarily in the
context of $v$-coherent domains, and we obtain quite satisfactory
results for Pr\"ufer $v$-multiplication domains.

Section 3 is devoted to a study of the $t\#$-property for
polynomial rings.  We show that if $R$ has the $t\#$-property,
then so does $\rx$ and that the converse is true if $\rx$ is
assumed to be $v$-coherent.  We also consider the $t\#$-property
in two commonly studied localizations of $\rx$.




\section{The $t\#$-property}

For a nonzero fractional ideal $I$ of a domain $R$ with quotient
field $K$, we set $I^{-1}=(R:_K I)=\{x \in K \mid xI \sub R\}$,
$I_v=(I^{-1})^{-1}$, and $I_t=\bigcup J_v$, where the union is
taken over all nonzero finitely generated subideals $J$ of $I$.
The reader is referred to \cite{g1} for properties of these (and
other) star operations.  We also recall that $I$ is said to be
\emph{divisorial} if $I=I_v$ and to be a \emph{$t$-ideal} if
$I=I_t$.  Finally, we denote the set of maximal $t$-ideals of $R$
by $\tmax(R)$.

We begin by repeating the definition of the $t\#$-property.

\begin{definition} A domain $R$ has the \emph{$t\#$-property} (or is a
\emph{$t\#$-domain}) if $\bigcap_{M \in \mc M_1} R_M \ne
\bigcap_{M \in \mc M_2} R_M$ for any two distinct subsets $\mc
M_1$ and $\mc M_2$ of $\tmax(R)$.
\end{definition}

\begin{theorem} \label{t:sharpchar} The following statements are
equivalent for a domain $R$.

\begin{itemize}

\item[(1)] $R$ is a $t\#$-domain.

\item[(2)] For each $N \in t\Max (R)$, we have $\bigcap_{M
\in t\Max (R) \sm \{N\}}R_M \nsubseteq R_N$.

\item[(3)] For each $N \in t\Max (R)$, we have $R \ne \bigcap_{M
\in t\Max (R) \sm \{N\}}R_M$.

\item[(4)] $\bigcap_{M \in \mc M_1} R_M$ and $\bigcap_{M \in \mc
M_2}R_M$ are incomparable for each pair of disjoint subsets $\mc
M_1$ and $\mc M_2$ of $t\Max (R)$.

\item[(5)] For each maximal $t$-ideal $M$ of $R$, there is a
divisorial ideal of $R$ which is contained in $M$ and no other
maximal $t$-ideal of $R$.

\item[(6)] For each maximal $t$-ideal $M$ of $R$, there is an
element $u \in K \sm R$ such that $M$ is the only maximal
$t$-ideal containing $(R:_R u)$.
\end{itemize}
\end{theorem}

\begin{proof}  By \cite[Proposition 4]{gr} we have $R=\bigcap_{M \in
t\text{-Max}(R)}R_M$.  Using this and the definitions, the
following implications are straightforward: (1) $\lra$ (2), (2)
$\lra$ (4); (2) $\lra$ (3); and (6) $\ra$ (5).

To prove (2) $\lra$ (6), observe that an element $u \in K$
satisfies $u \in \bigcap_{M \in t\text{-Max}(R) \sm \{N\}}R_M \sm
R_N$ if and only if $(R:_R u)$ is contained in $N$ and no other
maximal $t$-ideal of $R$.

Now assume (5).  Let $N \in t\text{-Max}(R)$, and pick a
divisorial ideal $I$ with $I \sub N$ and $I \nsubseteq M$ for each
$M \in t\text{-Max}(R) \sm \{N\}$.  Then $I^{-1} \sub \bigcap_{M
\in t\text{-Max}(R) \sm \{N\}}R_M$ (since for each $M$, we have
$(R:_R I^{-1}) = I \nsubseteq M$) , but $I^{-1} \nsubseteq R$.
Hence (5) $\ra$ (3), and the proof is complete.
\end{proof}

\begin{corollary} \label{c:tv} If $R$ is a domain with the
property that each maximal $t$-ideal is divisorial, then $R$ is a
$t\#$-domain.
\end{corollary}

\begin{proof} This is clear from the equivalence of conditions (1)
and (5) of Theorem~\ref{t:sharpchar}.
\end{proof}

Recall that a \emph{Mori domain} is a domain satisfying the
ascending chain condition on (integral) divisorial ideals.
Equivalently, a domain $R$ is a Mori domain if for each ideal $I$
of $R$ there is a finitely generated ideal $J \sub I$ with
$I_v=J_v$.  In particular, the $v$- and $t$-operations on a Mori
domain are the same. Hence Corollary~\ref{c:tv} implies that Mori
domains are $t\#$-domains.

\begin{remark} \label{r:sharpchar} We have stated
Theorem~\ref{t:sharpchar} for the $t$-operation, since our primary
interest is in that particular star operation.  However, suppose
that for a finite-type star operation $*$, we call a domain $R$ a
\emph{$*\#$-domain} if for each pair of nonempty subsets $\mc M_1$
and $\mc M_2$ of $*\text{-Max}(R)$ with $\mc M_1 \ne \mc M_2$, we
have $\bigcap_{M \in \mc M_1}R_M \ne \bigcap_{M \in \mc M_2}R_M$.
Then Theorem~\ref{t:sharpchar} remains true with $t$ replaced
everywhere (including the proof) by $*$.  This is of some interest
even in the case where $*$ is the trivial star operation ($I^*=I$
for each ideal $I$; this is often referred to as the
$d$-operation).  For the trivial star operation, Olberding has
proved the equivalence of statements (1) and (5) in
\cite[Proposition 2.2]{o}.
\end{remark}

Theorem~\ref{t:sharpchar} (1) $\lra$ (6) generalizes \cite[Theorem
1 (a) $\lra$ (b)]{gh}, which states that a Pr\"ufer domain $R$ is
a $\#$-domain if and only if each maximal ideal of $R$ contains a
finitely generated ideal which is contained in no other maximal
ideal of $R$.  This follows upon recalling that for $R$ Pr\"ufer
(i) each ideal is a $t$-ideal (so that
$t\text{-Max}(R)=\text{Max}(R)$) and (ii) for each $u \in K$,
$(R:_Ru)$ is finitely generated (in fact, two generated). In
general, one cannot hope to show that each maximal $t$-ideal of a
$t\#$-domain contains a finitely generated ideal which is
contained in no other maximal $t$-ideal, as the following example
shows.  (Example~\ref{e:ho} below is another such example.
However, that example has (Krull) dimension two, and we think it
might be of some interest to have a one-dimensional example.)

\begin{example} \label{e:pullback} Let $T$ be an almost Dedekind
domain with exactly one noninvertible maximal ideal $M$. (One such
example is constructed in \cite[Example 42.6]{g1}.)  By
\cite[Theorem 3]{g3}, $T$ is not a $\#$-domain.  For our purposes,
it does no harm to assume that $T/M$ has a proper subfield.  This
follows from the fact that $T(X)=T[X]_S$, where $S$ is the
multiplicatively closed subset of $T[X]$ consisting of those
polymomials $g$ having unit content (the ideal generated by the
coefficients of $g$), is also an almost Dedekind domain with
exactly one nonivertible maximal ideal, namely $MT(X)$, whose
residue field $T(X)/MT(X) \approx (T/M)(X)$ has infinitely many
proper subfields \cite[Proposition 36.7]{g1}. Let $F$ be such a
proper subfield of $T/M$, and let $R$ be defined by the following
pullback diagram of canonical homomorphisms: 

\[
\begin{CD}
		R   @>>>    F\\
		@VVV        @VVV    \\
        T  @>>>    k=T/M.\\
\end{CD} \]  

We claim that $R$ is a $t\#$-domain.  (In fact,
since $R$ is one dimensional, it is a $\#$-domain.)  We show that
$R$ satisfies condition (5) of Theorem~\ref{t:sharpchar}.  For
this it suffices to observe that each maximal ideal of $R$ is
divisorial.  This is clear for $M$, and if $P$ is a maximal ideal
of $R$ with $P \ne M$, then by \cite[Theorem 2.35]{gho} $P$ is
actually invertible.  Hence $R$ is a $(t)\#$-domain. Since $T$ is
a non-$t\#$-Pr\"ufer domain with offending maximal ideal $M$,
however, there is no finitely generated ideal of $T$ contained in
$M$ but no other maximal ideal of $T$; clearly, a similar
statement applies to $R$.
\end{example}

If we restrict our attention to domains in which conductors are
required to be finitely generated, i.e., to \emph{finite conductor
domains}, then the $t\#$-property becomes equivalent to the
property that each maximal $t$-ideal contain a finitely generated
ideal contained in no other maximal $t$-ideal.  In fact, we can
obtain such a result by requiring a little less than finite
generation of conductors.  Recall that a domain $R$ is said to be
\emph{$v$-coherent} if for each finitely generated ideal $I$ of
$R$, $I^{-1}$ has finite type (i.e., there is a finitely generated
ideal $J$ with $I^{-1}=J_v$). This condition was first studied
(under a different name) by Nour el Abidine \cite{n}.  It is easy
to see that a finite conductor domain is $v$-coherent.  We have
the following result.

\begin{theorem} \label{t:sharpchar2} For a $v$-coherent domain, the
conditions of Theorem~\ref{t:sharpchar} are each equivalent to:
Each maximal $t$-ideal of $R$ contains a finitely generated ideal
which is contained in no other maximal $t$-ideal of $R$.
\end{theorem}

\begin{proof}  The stated condition clearly implies condition (5)
of Theorem~\ref{t:sharpchar}.  On the other hand, condition (6),
in the presence of $v$-coherence, implies the stated condition.
\end{proof}

Now \cite[Theorem 1]{gh} contains a third equivalence, namely that
$R$ is uniquely representable as an intersection of a family
$\{V_{\al}\}$ of valuation overrings such that there are no
containment relations among the $V_{\al}$.  Since each valuation
overring of a Pr\"ufer domain is a localization, this suggests
exploring the possibility that the $t\#$-property on a domain $R$
is equivalent to the condition that $R$ contain a unique set of
incomparable $t$-primes $\{P_{\al}\}$ such that $R=\bigcap
R_{P_{\al}}$.  One implication is easy. If we assume the existence
of a unique set of $t$-primes $\{P_{\alpha}\}$ such that
$R=\bigcap R_{P_{\alpha}}$, then that set must be $\tmax(R)$, and
so $R$ is a $t\#$-domain. In Theorem~\ref{t:sharpchar3}, we
provide a converse in two cases.  First, we give an example
showing that the converse does not hold in general. Recall that a
domain $R$ is a \emph{Pr\"ufer $v$-multiplication domain} (PVMD)
if $R_M$ is a valuation domain for each maximal $t$-ideal $M$ of
$R$.

\begin{example} \label{e:ho} In \cite{ho} Heinzer and Ohm give an
example of an essential domain $D$ which is not a PVMD.  In their
example $k$ is a field, and $y,z,x_1, x_2, \ldots$ are independent
indeterminants over $k$; $R=k(x_1, x_2, \ldots)[y,z]_{(y,z)}$; for
each $i$, $V_i$ is a rank one discrete valuation ring containing
$k(\{x_j\}_{j \ne i})$ such that $y$, $z$, and $x_i$ all have
value 1; and $D=R \cap (\bigcap_i V_i)$. Then $D$ is a
2-dimensional domain, and in \cite{mz} it is shown that the
maximal ideals of $D$ are $M, P_1, P_2, \ldots$, where $M$ is the
contraction of the maximal ideal of $R=D_M$, and $P_i$ is the
contraction of the maximal ideal of $V_i=D_{P_i}$.  Note that each
$P_i$ has height one and is therefore a $t$-ideal. We observe that
each element of $R$ is also in $V_i$ for all but finitely many
$i$; this is the case since an element of $R$ involves only
finitely many of the $x_j$, and $x_i$ is a unit of $V_j$ for all
$j \ne i$. Similarly, each element of the maximal ideal of $R$ is
in the maximal ideal of $V_i$ for all but finitely many $i$.  It
follows that if $I$ is a finitely generated ideal of $D$ contained
in $M$, then $I$, and hence also $I_t$, is contained in all but
finitely many of the $P_i$.  Suppose that for such an $I$ we have
$I_t \nsubseteq M$. Write $1=x+m$ with $x \in I_t$ and $m \in M$.
By the observations stated above, $x$ and $m$ must be
simultaneously in all but finitely many of the $P_i$, a
contradiction.  Thus $M$ is a $t$-ideal.\footnote{We observe that
since $D_M=R$ is not a valuation domain, this yields an easy way
to see that $D$ is not a PVMD.} We show that $R$ has the
$t\#$-property by showing that $M$ and the $P_i$ satisfy condition
(5) of Theorem~\ref{t:sharpchar}. For each $i$, the divisorial
ideal $x_iD \subseteq P_i$, while $x_i \notin M$ and $x_i \notin
P_j$ for $j \ne i$. As for $M$, note that $y/z \in \bigcap D_{P_i}
\setminus D_M$.  Hence $(D:_D y/z) \sub M$, but $(D:_D y/z)
\nsubseteq P_i$ for $i=1,2, \ldots$. Finally, denoting the
height-one primes contained in $M$ by $\{Q_{\alpha}\}$, we observe
that $D=(\bigcap_{\al}D_{Q_{\al}}) \cap (\bigcap_i D_{P_i})$ for
the set of incomparable $t$-primes $\{Q_{\al}\} \cup \{P_i\} \ne
\tmax(R)$
\end{example}

In our next result, we use the fact that a PVMD is $v$-coherent
\cite{n}.

\begin{theorem} \label{t:sharpchar3}  Let $R$ be a either a
PVMD or a Mori domain. Then the conditions of
Theorem~\ref{t:sharpchar} are each equivalent to: There is a
unique set $\{P_{\al}\}$ of incomparable $t$-primes such that
$R=\bigcap R_{P_{\al}}$.  In particular, a Mori domain has this
property.
\end{theorem}

\begin{proof}  One implication was discussed above.  Assume
that $R$ is a $t\#$-domain, and suppose that $R=\bigcap
R_{P_{\al}}$ for some set $\{P_{\al}\}$ of incomparable
$t$-primes.  Observe that the hypotheses guarantee that $R$ is
$v$-coherent. It suffices to show that each $P_{\al}$ is a maximal
$t$-ideal. By way of contradiction, suppose that $P_{\be}
\subsetneqq M$, where $M$ is a maximal $t$-ideal. If $R$ is a
PVMD, then (since the $P_{\al}$ are incomparable and $R_M$ is a
valuation domain), $P_{\al} \nsubseteq M$ for each $\al \ne \be$.
By Theorem~\ref{t:sharpchar2}, $M$ contains a finitely generated
ideal $I$ which is contained in $M$ and no other maximal
$t$-ideal.  In particular, $I \nsubseteq P_{\al}$ for $\al \ne
\be$. Pick $a \in M \sm P_{\be}$.  Then $(I,a)$ is a finitely
generated ideal contained in no $P_{\al}$ whatsoever.  It follows
that $(I,a)^{-1} \sub \bigcap R_{P_{\al}} = R$, whence
$(I,a)_v=R$. However, since $M$ is a $t$-ideal, we have $(I,a)_v
\sub M$, a contradiction in this case.  If $R$ is Mori, then $M$
itself is divisorial, and, since $M$ is contained in no $P_{\al}$,
we obtain the contradiction that $M^{-1} \sub \bigcap
R_{P_{\al}}$.
\end{proof}

\begin{remark} \label{r:vcohsharpchar} We have not been able to
determine whether weakening the hypothesis of
Theorem~\ref{t:sharpchar3} to $v$-coherent is sufficient.  It does
suffice if we make the following subtle change to the condition:
there is a unique set $\{P_{\al}\}$ of $t$-primes such that both
$R=\bigcap R_{P_{\al}}$ and the intersection is irredundant (no
$R_{P_{\al}}$ can be deleted). To see this, suppose that $R$ is a
$t\#$-domain, and let $\{P_{\al}\}$ be as indicated. Pick a
$P_{\be}$; we wish to show that it is a maximal $t$-ideal. The
irredundancy hypothesis allows us to choose $u \in R_{P_{\be}} \sm
\bigcap_{\al \ne \be} R_{P_{\al}}$.  We have $(R:_R u) \nsubseteq
P_{\be}$ and $(R:_R u) \sub P_{\al}$ for each $\al \ne \be$. Since
$R$ is $v$-coherent, there is a finitely generated ideal $I$ with
$(R:_R u)=I_v$. Pick a maximal $t$-ideal $M \supseteq P_{\be}$. If
there is an element $a \in M \sm P_{\be}$, then, as in the proo
of Theorem~\ref{t:sharpchar3}, the ideal $(I,a)$ will furnish a
contradiction.
\end{remark}




\section{Overrings of $t\#$-domains}

In \cite{gh} Gilmer and Heinzer also studied Pr\"ufer domains with
the property that each overring is a $\#$-domain; these domains
have come to be called \emph{$\#\#$-domains}.  Our goal in this
section is to obtain $t$-analogues of results on $\#\#$-domains.

Most of the characterizations of Pr\"ufer $\#\#$-domains in
\cite{gh} can be extended to PVMDs with the property that each
$t$-linked overring is $t\#$. However, if we want to consider a
larger class of domains, e.g., $v$-coherent domains, the question
arises as to which overrings should be considered. Put another
way, it is not clear exactly how one should define the
$t\#\#$-property (and we shall not do so).

In what follows, it will be convenient to employ the language of
localizing systems. We recall the requisite definitions.  A
nonempty set $\mc F$ of nonzero ideals of $R$ is said to be a
\emph{multiplicative system of ideals} if $IJ \in \mc F$ for each
$I,J \in \mc F$.  The ring $R_{\mc F}= \{x \in K \mid xI \sub R
\text{ for some } I \in \mc F\}$ is called a \emph{generalized
ring of quotients} of $R$. For each ideal $J$ of $R$ we set
$J_{\mc F}=\{x \in K \mid xI \sub J \text{ for some }I \in \mc
F\}$; $J_{\mc F}$ is an ideal of $R_{\mc F}$ containing $JR_{\mc
F}$.

A particular type of multiplicative system is a \emph{localizing
system}: this is a set $\mc F$ of ideals of $R$ such that (1) if
$I \in \mc F$ and $J$ is an ideal of $R$ with $I \sub J$, then $J
\in \mc F$ and (2) if $I \in \mc F$ and $J$ is an ideal of $R$
such that $(J:_R a) \in \mc F$ for every $a \in I$, then $J \in
\mc F$.  If $\La$ is a subset of $\op{Spec} R$, then $\mc
F(\La)=\{I \mid \text{ is an ideal of $R$ such that }I \nsubseteq
P \text{ for each } P \in \La\}$ is a localizing system; moreover,
$R_{\mc F(\La)}=\bigcap_{P \in \La}R_P$. A localizing system $\mc
F$ is said to be \emph{spectral} if $\mc F = \mc F(\La)$ for some
set of primes $\La$.  Finally, an \emph{irredundant spectral
localizing system} is a localizing system of ideals $\mc F(\La)$,
where $\La$ is a set of pairwise incomparable primes.

These notions have $t$-analogues.  A set of $t$-ideals is a
\emph{$t$-multiplicative system} if it is closed under
$t$-multiplication; a $t$-multiplicative system $\Phi$ is a
\emph{$t$-localizing system} if it satisfies the closure
operations (1) and (2) above.

 The localizing system
$\mc F$ is said to be of \emph{finite type} if for each $I \in \mc
F$ there is a finitely generated ideal $J \in \mc F$ with $J \sub
I$. Also, $\mc F$ is said to be \emph{$v$-finite} if each
$t$-ideal of $\mc F$ contains a $v$-finite ideal which is also in
$\mc F$.

Denoting the set of $t$-ideals of $R$ by $t(R)$, it is easy to see
that if $\mc F$ is a  localizing system, then $\Phi=\mc F \cap
t(R)$ is a $t$-localizing system, $R_{\mc F}=R_{\Phi}$, and $\mc
F$ is $v$-finite if and only if $\Phi$ is $v$-finite.  Conversely,
if $\Phi$ is a $t$-localizing system of $t$-ideals, then
$\ol{\Phi}=\{I \mid I_t \in \Phi\}$ is a localizing system of
ideals with $\Phi=\ol{\Phi} \cap t(R)$.

Let $\La$ be a set of pairwise incomparable $t$-primes.  With $\mc
F(\La)$ as above and  $\Phi(\La)=\mc F(\La)\cap t(R) \,(\, = \{I
\mid I \text{ is a t-ideal and } I \nsubseteq P \text{ for all } P
\in \La\})$, we have that   $I\in \mc F(\La)$ if and only if
$I_{t} \in \Phi(\La)$. Hence $\overline{\Phi(\La)}=\mc F(\La)$.

An overring $T$ of $R$ is a \emph{$t$-subintersection of $R$} if
it has the form $\bigcap R_P$, where the intersection is taken
over some set of $t$-primes $P$ of $R$, i.e., if $T =
R_{\Phi(\La)}$ for some spectral $t$-localizing system $\Phi(\La)$
of $R$, where $\La$ is a set of $t$-primes. We say that $T$ is
\emph{$t$-flat over $R$} if $T_M=R_{M \cap R}$ for each maximal
$t$-ideal $M$ of $T$ \cite{kp}. Finally, recall that $T$ is
$t$-linked over $R$ if for each finitely generated ideal $I$ of
$R$ with $(R:I)=R$ we have $(T:IT)=T$ \cite{dhlz}.

The following implications are easily verified: $T$ is $t$-flat
over $R$ $\ra$ $T$ is a $t$-subintersection of $R$ $\ra$ $T$ is a
generalized ring of quotients of $R$ $\ra$ $T$ is $t$-linked over
$R$.

All these conditions are equivalent for PVMDs \cite[Proposition
2.10]{kp}, but we believe that in general none of the arrows can
be reversed if $R$ is merely assumed to be $v$-coherent. Also, if
$R$ is a PVMD, then every $t$-linked overring of $R$ is a PVMD
\cite[Theorem 3.8 and Corollary 3.9]{k}, but if $R$ is just
$v$-coherent, we know only that generalized rings of quotients of
$R$ are $v$-coherent \cite[Proposition 3.1]{gab2}.

We shall begin by considering $t$-flat overrings of $v$-coherent
domains. Recall that, for any domain $R$, an overring $T$ of $R$
is $t$-flat over $R$ if and only if $T$ is a generalized ring of quotients
with respect to a $v$-finite $t$-localizing system of ideals
\cite[Theorem 2.6]{elb1}.

On the other hand, we know that if $R$ is  Pr\"ufer then every
overring is flat, and we also know that $R$ is a $\#\#$-domain iff
each irredundant spectral localizing system is finitely generated
\cite{gab1}.  We shall show that for $v$-coherent domains the
property that each irredundant spectral $t$-localizing system is
$v$-finite is equivalent to the property that each
$t$-subintersection of $R$ is $t$-flat and $t\#$.

\begin{lemma}\label{lem}  Let $R$ be a $v$-coherent domain and $\Phi$  a
    $t$-localizing system of $t$-ideals.
    Then the following statements are equivalent.
    \begin{itemize}
\item[(1)] $\Phi$  is $v$-finite.
\item[(2)] The set $\La$ of maximal elements of $t\Spec (R)
\setminus \Phi$ is not empty, and $M \in t\Max (R_{\Phi})$
if and only if $M = P_{\Phi}$ for some $P\in \La$.
\end{itemize}
Under these conditions, $\Phi=\Phi(\La)$. In particular, if $\La$
is a set
    of pairwise incomparable $t$-primes of $R$, then $\Phi=\Phi(\La)$  is $v$-finite
    if and only if  $t\Max (R_{\Phi}) =\{P_{\Phi} \mid P \in \La\}$.
\end{lemma}

\begin{proof}
Set $\mc F = \ol \Phi= \{I \mid I_{t}\in \Phi\}$ and use (i)
$\lra$ (vi) of \cite[Theorem 3.3]{gab2}.
 \end{proof}

\begin{proposition}\label{1} Let $R$ be a $v$-coherent domain. If $\La$ and $\La'$
        are two sets of pairwise incomparable $t$-primes such that
        $\Phi(\La)$ and $\Phi(\La')$ are $v$-finite and $R_{\Phi(\La)}
       = R_{\Phi(\La')}$, then $\La = \La'$.
	\end{proposition}

    \begin{proof} By Lemma~\ref{lem}, we have $\tmax(T) =\{P_{\Phi} \mid P\in\La\} =
        \{Q_{\Phi'} \mid Q\in \La'\}$ and, upon contracting to $R$, we
 obtain $\La = \La'$.
        \end{proof}

            Recalling that an overring $T$ of a domain $R$ is $t$-flat over $R$
            if and only if
     $T=R_{\Phi}$ for some  $v$-finite $t$-localizing
    system $\Phi$, the preceding two results immediately imply:

    \begin{corollary}\label{cor} Let $R$ be a $v$-coherent domain and let $T$ be a
    $t$-flat overring of $R$. Then there exists a uniquely determined set $\La$
    of pairwise incomparable $t$-primes for which $T=R_{\Phi(\La)}$ and
    $\Phi(\La)$ is $v$-finite. The set $\La$ is given by
    $\La = \{M\cap R \mid M \in t\Max (T)\}$.
    \end{corollary}

    \begin{proposition}\label{2}  Let $R$ be a $v$-coherent domain.
    Then the following statements are equivalent.
\begin{itemize}

\item[(1)] For each set $\La$
of pairwise incomparable $t$-primes of $R$, $\Phi(\La)$ is $v$-finite.

\item[(2)] If $\La$ and $\La'$
        are two sets of pairwise incomparable $t$-primes of $R$ such that
     $R_{\Phi(\La)} = R_{\Phi(\La')}$, then $\La = \La'$.

    \item[(3)] If $T$ is a $t$-subintersection of $R$ and is
    represented as $T=\bigcap_{P \in \La}R_P$ for some set $\La$ of
    pairwise incomparable $t$-primes, then that representation is irredundant.

    \item[(4)]  For each $t$-prime $P$ and each set $\La$
    of pairwise incomparable $t$-primes of $R$ not containing $P$,
there exists an element $u\in K$ such that $(R:_{R}u)\subseteq P$
and $(R:_{R}u)\nsubseteq Q$, for each $Q\in\La$.

    \item[(5)]  For each $t$-prime $P$ and each set $\La$
    of pairwise incomparable $t$-primes of $R$ not containing $P$,
there exists a finitely generated ideal $J$ of $R$ such that $J
\sub P$ and $J \nsubseteq Q$ for each $Q \in \La$.

    \item[(6)] For each $t$-prime $P$ and each set $\La$
    of pairwise incomparable $t$-primes of $R$ not containing $P$,
$R_{P}\nsupseteq R_{\Phi(\La)}$.

\item[(7)] For each set $\La$
    of pairwise incomparable $t$-primes of $R$, $R_{\Phi(\La)}$
    is $t$-flat over $R$ and has the $t\#$-property.
\end{itemize}
\end{proposition}

\begin{proof}
    (1) $\ra$(2) by Proposition \ref {1}.

	(2) $\ra$ (3) is clear.

	(3) $\ra$ (1): Given a set $\La$ of incomparable primes, consider
the $t$-subintersection $T=R_{\Phi(\La)}$ of $R$. Since $T$ is
$v$-coherent \cite[Proposition 3.1]{gab2} and the intersection is
irredundant, we obtain $t$-Max$(T) =\{P_{\Phi(\La)} \mid P \in
\La\}$ as in Remark~\ref{r:vcohsharpchar}. It follows that
$\Phi(\La)$ is $v$-finite (Lemma~\ref{lem}).

    (1) $\ra$ (7):   Let $\La$ be a set
    of pairwise incomparable $t$-primes of $R$ and $T= R_{\Phi(\La)}$.
Since $\Phi(\La)$ is $v$-finite, then $T$ is $t$-flat over $R$,
and $\La = \{M\cap R \mid M \in \tmax(R)\}$ is uniquely determined
by Corollary \ref{cor}. Hence we cannot delete any $P\in \La$, and
so the intersection is irredundant. In addition, by $t$-flatness,
$T_{M}=R_{M\cap R}$; hence $T$ is a $t\#$-domain.

    (7) $\ra$ (3): If the $t$-subintersection
    $T= R_{\Phi(\La)}$ of $R$ is $t$-flat, then
    $\La = \{M\cap R \mid M \in t$-Max$(R)\}$ by Corollary \ref{cor}. If $T$
    is also a $t\#$-domain, then
     $T = \bigcap T_{M}=\bigcap R_{M\cap R}$ is
    an irredundant $t$-subintersection.

    (2) $\ra$ (4): Given $\La$ and $P$ as specified, set
    $\La'=(\La \sm \{Q \in \La \mid Q \sub P\}) \cup \{P\}$.  Then
    $\La \ne \La'$, so that $R_{\Phi(\La)} \ne R_{\Phi(\La')}$ by
    (2).  Since we clearly have $R_{\Phi(\La')} \sub
    R_{\Phi(\La)}$, there is an element $u \in R_{\Phi(\La)} \sm
    R_{\Phi(\La')}$, and for this $u$ we have $(R:_R u) \sub P$ and
    $(R:_{R} u) \nsubseteq Q$ for each $Q \in \La$.

    (4) $\ra$ (5): Since $R$ is $v$-coherent, then the ideal $(R:_{R}u)$
    contains a finitely generated subideal $J$ with $J_v=(R:_{R}u)$; this
    $J$ does what is required.

    (5) $\ra$ (6): Given $J$ as indicated, one shows easily that
    $(R:J) \sub R_{\Phi(\La)}$ but $(R:J) \nsubseteq R_P$, whence
    $R_P \nsupseteq R_{\Phi(\La)}$.

    (6) $\ra$ (2):  Suppose that $\La$ and $\La'$ are two sets of
    pairwise incomparable primes for which
    $R_{\Phi(\La)}=R_{\Phi(\La')}$ but $\La \ne \La'$.  We may
    then assume that there is a prime $P \in \La \sm \La'$.  If $P
    \nsubseteq Q$ for all $Q \in \La'$, then (6) yields $R_P
    \nsubseteq R_{\Phi(\La')}=R_{\Phi(\La)}$, a contradiction.
    We then denote by $\La''$ the maximal elements in the set
    $(\La \cup \{Q \in \La' \mid P \sub Q\}) \sm \{P\}$.  (Choosing
    the maximal elements is possible since both $\La$ and $\La'$ contain
    pairwise incomparable elements.) Pick
    $Q_0 \in \La'$ with $P \sub Q_0$.  Then $Q_0 \in \La''$, and
    we have $R_P \supseteq R_Q \supseteq R_{\Phi(\La'')}$, which
    contradicts (6).
\end{proof}

\begin{remark} \label{r:mori} The equivalent conditions of
Proposition~\ref{2} hold automatically for a Mori domain, since in
such a domain each $t$-ideal is $v$-finite.
\end{remark}

\begin{proposition}\label{2bis}  Let $R$ be a $v$-coherent domain such that
    each $t$-subintersection of $R$ is $t$-flat over $R$.
Then the following statements are equivalent.
\begin{itemize}

    \item[(1)]  Each $t$-flat overring of $R$ is a $t\#$-domain.

    \item[(2)]  Each  $t$-subintersection of $R$ is a $t\#$-domain.

\item[(3)]  For each set $\La$
    of pairwise incomparable $t$-primes of $R$, $\Phi(\La)$ is $v$-finite.

 \item[(4)] If $T$ is a $t$-subintersection of $R$, there exixts a unique
set of pairwise incomparable $t$-primes $\La$ of $R$ such that $T
= R_{\Phi(\La)}$; moreover, $\La = \{M\cap R \mid M \in t\Max
(T)\}$.

\item[(5)] If $T$ is a $t$-flat overring of $R$ and
$T= \bigcap_{Q \in \Lambda} T_Q$ for some set $\Lambda$ of
pairwise incomparable $t$-primes of $T$, then
$\Lambda=t\Max (T)$.
\end{itemize}
\end{proposition}

\begin{proof} (1)$\lra$(2) and (5)$\ra$(1) are clear.

    (2)$\lra$(3) by Proposition \ref{2}.

    (3)$\ra$(4) by Corollary \ref{cor}.

    (4) $\ra$ (5): Assume that $T$ is a $t$-flat overring of $R$
    and that we have $T=\bigcap_{Q \in \Lambda} T_Q$, where $\Lambda$ is a set
    of pairwise incomparable $t$-primes of $T$.  By $t$-flatness,
    $T_Q = R_{Q \cap R}$ for each $Q \in \La$.  Hence
    $T=R_{\Phi(\Gamma)}$, where $\Gamma =\{Q \cap R \mid Q \in
    \Lambda\}$.  We then have $\Lambda=\tmax(T)$ by (4) (and
    $t$-flatness).
    \end{proof}

If $R$ is a Mori domain, then, as mentioned in
Remark~\ref{r:mori}, the equivalent conditions of
Proposition~\ref{2} hold.  It then follows from \cite[Theorem
2.6]{elb1} that each $t$-subintersection of $R$ is $t$-flat; hence
the equivalent conditions of Proposition~\ref{2bis} hold also.

For a PVMD, t-subintersections are automatically t-flat; in fact,
t-linked overrings are t-flat by \cite[Proposition 2.10]{kp}. Thus
the hypotheses of Proposition~\ref{2bis} hold for PVMDs. Our next
proposition adds several more equivalences for PVMDs.  We need the
following lemma.

\begin{lemma}\label{dual} Let $R$ be a PVMD and let $P$ be a $t$-prime of
$R$ which is not $t$-invertible. Then $(P:P) = (R:P) = R_P \cap
S$, where $S = \bigcap_{M \in t\Max (R), M \nsupseteq P}
R_{M}$.
\end{lemma}

\begin{proof}
    By \cite[Proposition 2.3 and Lemma 1.2]{hz}, $(R:P)=(P:P)$.  The
result now follows from \cite [Theorem 4.5]{hklm}.
\end{proof}

\begin{proposition}\label{3} For a PVMD $R$, the conditions of
Proposition~\ref{2bis} are also equivalent to each of the
following.
\begin{itemize}

	\item[(6)]  If $\La \subseteq t\Max (R)$, then $\Phi(\La)$ is
$v$-finite.

    \item[(7)]  Each  $t$-prime ideal $P$ of $R$ contains a finitely generated ideal
which is not contained in any maximal $t$-ideal of $R$ not
containing $P$.

\item[(8)]  For each $t$-prime $P$ of $R$,
there exists an element $u\in K$ such that $(R:_{R}u)\subseteq P$
and $(R:_{R}u)\nsubseteq M$, for each maximal $t$-ideal $M$ not
containing $P$.

\item[(9)] For each $t$-prime ideal $P$ of $R$, $R_{P}\nsupseteq
\bigcap R_{M}$, where $M$ ranges over the set of maximal
$t$-ideals not containing $P$.

\item[(10)] Each $t$-linked overring of $R$ is a $t\#$-domain.

\item[(11)] $(P:P)$ is a $t\#$-domain for each  $t$-prime $P$ of $R$.

\end{itemize}
\end{proposition}

\begin{proof}
(3) $\ra$ (6) is clear.

    (6) $\ra$ (7): If $\La$ is the set of maximal $t$-ideals not containing
    $P$, then $P\in \Phi(\La)$ and $\Phi(\La)$ is $v$-finite.

    (7) $\ra$ (1): Let $T$ be a $t$-subintersection of $R$.
Then $T$ is $t$-flat over $R$, and we have $T =
    \bigcap_{M \in \tmax(T)}R_{M\cap R}$. Fix $N \in t$-Max($T$)
    and let $J$ be a finitely generated ideal of $R$ contained in $P=N\cap R$ and
    not contained in the maximal $t$-ideals of $R$ not containing $P$.
    Since in a PVMD two incomparable $t$-primes are $t$-comaximal,
    then $J$ is not contained in $M\cap R$ for each maximal $t$-ideal
    $M\neq N$ of $T$. It follows that $JT$ is a finitely generated
    ideal contained in $N$ and not contained in $M$ for $M\neq N$. We
    conclude by applying Theorem~\ref{t:sharpchar2}.

(3) $\ra$ (8) by Proposition \ref{2}.

(8) $\ra$ (7) by $v$-coherence.

(8) $\lra$ (9) because, for each prime $P$ and $u \in K$,
       $(R:_{R}u)\subseteq P$ iff $u \notin R_{P}$.

(1) $\lra$ (10) because each $t$-linked overring of a PVMD is
$t$-flat \cite[Proposition 2.10]{kp}.

(11) $\ra$ (9): Let $T = (P:P)$. If $P$ is $t$-invertible then $R
= T$.  Otherwise,  $T = (R:P) = R_{P} \cap (\cap R_{M_{\alpha}})$,
where $M_\alpha$ ranges over the set of maximal $t$-ideals of $R$
not containing $P$ (Lemma \ref{dual}). In either case, setting
$\La = \{P\} \cup \{M_{\alpha}\}$, we have that $T=
R_{\Phi(\La)}$. Since $R$ is $v$-coherent, the set of ideals
$\{Q_{\Phi(\La)} = QR_{Q}\cap T; Q \in \La \}$ is a set of
incomparable $t$-primes of $T$ \cite[Proposition 3.2]{gab2}. For
each $Q \in \La$, we have $R_Q = T_{Q_{\Phi(\La)}}$ and by
hypothesis $T$ is a $t\#$-domain. Hence by Theorem 1.6
$R_{\Phi(\La)}$ is an irredundant intersection. It follows that
$R_P \nsupseteq \cap R_{M_{\alpha}}$.

(10) $\ra$ (11): According to \cite[Proposition 2.2 (5)]{dhlz},
$(A_v:A_v)$ is $t$-linked over $R$ for each ideal $A$ of $R$. In
fact, it is easy to see that replacing ``$v$'' by ``$t$'' in the
proof of that result shows that $(A_t:A_t)$ is $t$-linked.  In
particular, if $P$ is a $t$-prime of $R$, then $(P:P)$ is
$t$-linked.
\end{proof}

Comparing conditions (3) and (6) of Propositions~\ref{2bis} and
\ref{3}, we observe that for PVMDs one has to consider only
subsets of $\tmax(R)$ rather than all sets of incomparable
$t$-primes.

The equivalence of conditions (7) and (8) above is also proved in
 \cite[Lemma 3.6]{elb2}.
The equivalence of conditions (10) and (11)
 for Pr\"ufer domains is \cite[Proposition 2.5]{o}.

\medskip

When $R$ is Pr\"ufer, Proposition \ref{3} recovers \cite[Theorem
2.4]{gab1}. In \cite[Theorem 2.5]{gab1} it is also proved that for
Pr\"ufer domains the $\#\#$-condition is equivalent to the
$\#_{P}$-condition introduced by N. Popescu in \cite{p}. We recall
that $R$ is a $\#_{P}$-domain if, given two sets of prime ideals
$\La_{1} \ne \La_{2}$ with the property that $P+Q = R$ for each
pair of distinct ideals $P\in \La_{1}$ and $Q\in \La_{2}$, we have
$R_{\Phi(\La_{1})} \neq R_{\Phi(\La_{2})}$.

We can define the $t\#_{P}$-property analogously: $R$ is a
\emph{$t\#_{P}$-domain} if, given two sets of prime $t$-ideals
$\La_{1} \ne \La_{2}$ with the property that $(P+Q)_{t}=R$ for
each pair of distinct ideals $P\in \La_{1}$ and $Q\in \La_{2}$, we
have $R_{\Phi(\La_{1})} \neq R_{\Phi(\La_{2})}$.

We will show that, with this definition, \cite[Theorem 2.5]{gab1}
can be extended to PVMDs. Recall that, if $R$ is a PVMD, then for
any two incomparable prime $t$-ideals $P$ and $Q$ we have
$(P+Q)_{t}=R$ (since $R_M$ is a valuation domain for each maximal
$t$-ideal $M$ of $R$).

\begin{proposition} \label{p:irred} Let $R$ be a $v$-coherent domain,
and assume that the equivalent conditions of Proposition \ref{2}
are satisfied. Then $R$ is a $t\#_{P}$-domain.
\end{proposition}

\begin{proof} Let $\La_{1} \ne \La_{2}$ be two sets of prime
$t$-ideals of $R$ with the property that $(P+Q)_{t}=R$ for each
pair of distinct ideals $P\in \La_{1}$ and $Q\in \La_{2}$, and let
$P \in \La_1 \setminus \La_2$.  Since $(P+Q)_t=R$ for $Q \in
\La_2$, we have $(P+M)_t=R$ for each $M$ in the set $\Ga=\{N \in
\tmax(R) \mid Q \sub N \text{ for some } Q \in \La_2\}$.  Since
$\Ga$ is a set of incomparable $t$-primes not containing $P$, we
may apply Proposition~\ref{2} (4) to obtain an element $u \in K$
such that $(R:_R u) \sub P$ but $(R:_R u) \nsubseteq M$ for each
$M \in \Ga$. It is then easy to see that $u \in R_{\Phi(\La_2)}
\setminus R_{\Phi(\La_1)}$.
\end{proof}

Our next result shows that for PVMDs the $t\#_P$-condition is
equivalent to the conditions of Propositions~\ref{2bis} and
\ref{3}.

\begin{proposition} \label{p:tsharpp} Let $R$ be a PVMD.
Then $R$ is a $t\#_{P}$-domain if and only if each $t$-linked
overring of $R$ is a $t\#$-domain.
\end{proposition}

\begin{proof} In a PVMD any two incomparable $t$-primes are
$t$-comaximal.  Hence if $R$ is a $t\#_P$-domain, then $R$ must
satisfy condition (3) of Proposition~\ref{2}.  The fact that
conditions (4) of  Proposition~\ref{2bis} and (10) of
Proposition~\ref{3} are equivalent then shows that each $t$-linked
overring of $R$ is a $t\#$-domain. The converse follows from
Proposition~\ref{p:irred}.
\end{proof}

The next result generalizes \cite[Theorem 2.6]{gab1}.

\begin{proposition}\label{acc} The following statements are
equivalent for a $v$-coherent domain $R$.
\begin{itemize}

    \item[(1)] For each set $\La$ of $t$-primes of $R$,
    $\Phi(\La)$ is $v$-finite.

\item[(2)] $R$ satisfies the ascending chain conditions on
$t$-primes, and $R$ satisfies the equivalent conditions of
Proposition~\ref{2}.
    \end{itemize}
    \end{proposition}

\begin{proof} (1) $\ra$ (2): Let $\La$ be a nonempty set of
$t$-primes of $R$.  Since $\Phi(\La)$ is $v$-finite,
Lemma~\ref{lem} implies that $\La$ has maximal elements.  Hence
$R$ satisfies the acc on $t$-primes.  Condition (1) of
Proposition~\ref{2} holds by hypothesis.

(2) $\ra$ (1): Let $\La$ be a nonempty set of $t$-primes.  Then
acc on $t$-primes implies that each element of $\La$ is contained
in a maximal element.  Hence if $\La_0$ is the set of maximal
elements of $\La$, then $\Phi(\La) = \Phi(\La_0)$ is $v$-finite by
Proposition~\ref{2}.
\end{proof}

The preceding result can be improved for PVMD's in a way which
generalizes \cite[Theorem 4]{gh}. We first recall some results
from \cite{elb1} and prove a variation on \cite[Lemma 4]{gh}.

\begin{lemma}\label{Elbag} Let $R$ be any domain. Then $R$ satisfies
the ascending chain condition on radical $t$-ideals if and only if
each prime $t$-ideal is the radical of a $v$-finite $t$-ideal.

If $R$ does satisfy the acc on radical $t$-ideals, then every
$t$-ideal has only finitely many minimal ($t$-)primes.
\end{lemma}

\begin{proof} \cite[Lemmas 3.7 and 3.8]{elb1}.
 \end{proof}

\begin{lemma}[{cf. \cite[Lemma 4]{gh}}] \label{l:gh}
Let $I=(a_1, \ldots, a_n)$ be a finitely generated ideal of a PVMD
$R$.  Then each minimal prime ideal of $I_v$ is minimal over some
$(a_i)$.  Moreover, if $I_v$ has only finitely many minimal
primes, then each minimal prime of $I_v$ is the radical of a
$v$-finite divisorial ideal.
\end{lemma}

\begin{proof} Let $P$ be minimal over $I_v$.  Then $P$ is a
$t$-prime, and, since primes contained in $P$ are also $t$-primes,
$P$ is also minimal over $I$.  The proof of the first statement
now proceeds as in the proof of the corresponding part of
\cite[Lemma 4]{gh}.  Now assume that $I_v$ has only finitely many
minimal primes $P_1, \ldots, P_k$, $k \ge 2$.  Since there are no
containment relations among the $P_i$ (and since the $t$-spectrum
of a PVMD is treed), we have $(P_1 + P_2 \cdots P_k)_t = R$. Hence
there are finitely generated ideals $A \sub P_1$ and $B \sub P_2
\cdots P_k$ with $(A+B)_v=R$.  We claim that $P_1$ is the radical
of $(I+A)_v$.  To see this, suppose that $Q$ is a prime which is
minimal over $(I+A)_v$.  Then $Q$ is a $t$-prime and must contain
a prime minimal over $I_v$; that is, $Q$ must contain one of the
$P_i$.  However, $Q$ cannot contain $P_i$ for $i \ge 2$, since
then $Q$ would contain $B$ (and $(A+B)_v=R$).  Hence $Q$ contains,
and is therefore equal to, $P_1$.
\end{proof}

\begin{proposition}\label{acc2} Let $R$ be PVMD. Then the statements in
    Proposition~\ref{acc} are equivalent to each of the following.
\begin{itemize}

    \item[(3)] $R$ satisfies
the ascending chain condition on radical $t$-ideals.

\item[(4)] $R$ satisfies the ascending chain condition on $t$-primes,
and, for each finitely generated ideal $I$, the set of minimal
primes of $I_v$ is a finite set.

    \item[(5)] Each $t$-prime of $R$ is branched and each $t$-linked
    overring of $R$ is a $t\#$-domain.

\item[(6)] $R$ satisfies the ascending chain condition on
$t$-primes and each $t$-linked overring of $R$ is a $t\#$-domain.
\end{itemize}
    \end{proposition}

    \begin{proof} (2) $\ra$ (4):
By (3) $\lra$ (5) of Proposition \ref{2},
        for each $t$-prime $P$ of $R$, we have that $R_{P}\nsupseteq
\bigcap R_{M}$, where the intersection is taken over those maximal
$t$-ideals of $R$ which do not contain $P$. Hence each principal
ideal has only finitely many minimal ($t$-)primes by \cite[Lemma
3.9]{elb2}.  Thus if $I=(a_1, \ldots, a_n)$ is finitely generated,
then $I_v$ can have only finitely minimal primes, since
Lemma~\ref{l:gh} implies that each such minimal prime must be
minimal over one of the $a_i$.

        (4) $\ra$ (3): Let $P$ be a $t$-prime of $R$.
By Lemma~\ref{Elbag}, it suffices to show that $P$ is the radical
of a $v$-finite $t$-ideal.  By the ascending chain condition on
$t$-primes, the set of $t$-primes properly contained in $P$ has a
maximal element $Q$. Thus, for $x \in P \setminus Q$, $P$ is
minimal over the principal ideal $xR$.  By assumption, $xR$ has
only finitely many minimal primes.  Hence Lemma~\ref{l:gh} yields that
$P$ is the radical of a $v$-finite $t$-ideal, as desired.

(3) $\ra$ (2): Clearly, $R$ satisfies the ascending chain
condition on $t$-primes. Let $P$ be a $t$-prime of $R$. By
Lemma~\ref{Elbag} $P$ is the radical of $J_{v}$ for some finitely
generated ideal $J$ of $R$. Since any $t$-prime containing $J$
also contains $P$, it is clear that condition (5) of
Proposition~\ref{2} holds.

(2) $\lra$ (6) because each $t$-linked overring of a PVMD is a
$t$-flat $t$-subintersection (\cite[Theorem 3.8]{k} and
\cite[Proposition 2.10]{kp}.

(5) $\lra$ (6): Since each localization of a PVMD at a $t$-prime
is a valuation domain, each $t$-prime of $R$ is branched if and
only if $R$ satisfies the ascending chain condition on $t$-primes.
    \end{proof}

    The PVMD's with the property that each $t$-localizing system of ideals
is $v$-finite have been studied in \cite{elb1} and \cite{elb2}.
They are called Generalized Krull domains. By \cite[Theorem
3.9]{elb1}, $R$ is a Generalized Krull domain if and only if each
principal ideal has only finitely many minimal primes and  $P \neq
(P^2)_{t}$ for each $t$-prime $P$. On the other hand, the first
condition is satisfied under the equivalent conditions of
Proposition \ref{3} \cite[Lemma 3.9]{elb2}. Hence we obtain the
following result.

\begin{corollary} A PVMD $R$ is a Generalized Krull domain if and only if each
$t$-linked overring of $R$ is a $t\#$-domain and $P \neq
(P^2)_{t}$ for each $t$-prime $P$.
\end{corollary}




\section{Polynomial rings}

In this section, we denote by $\{X_{\alpha}\}$ a set of
independent indeterminates over $R$.  Let us call a prime ideal
$Q$ of $\rx$ an \emph{upper to zero} if $Q \cap R=0$. For $f$ in
the quotient field of $\rx$, the \emph{content of $f$}, written
$c(f)$ is the fractional $R$-ideal generated by the coefficients
of $f$; we also write $c(I)$ for the fractional ideal generated by
the coefficients of all the polynomials in the fractional
$\rx$-ideal $I$.

\begin{lemma} \label{l:uppers} Let $Q$ be an upper to zero in $\rx$.
Then the following statements are equivalent.
\begin{itemize}
\item[(1)] $Q=fK[\{X_{\al}\}] \cap \rx$ for some irreducible polynomial $f \in
K[\{X_{\al}\}]$. (Note that we may take $f \in \rx$.)
\item[(2)] $\op{ht} Q=1$.
\item[(3)] $\rx_Q$ is a DVR.
\end{itemize}
\end{lemma}

\begin{proof} A localization argument establishes (1) $\lra$ (2),
and (3) $\ra$ (2) is trivial.  Assume (2).  If $\{X_{\al}\}$ is
finite, then (3) follows from a standard induction argument.  If
$\{X_{\al}\}$ is infinite, then we may pick $X_1, \ldots, X_n \in
\{X_{\al}\}$ with $Q \cap R[X_1, \ldots, X_n] \ne 0$.  Then
$V=R[X_1, \ldots, X_n]_{Q \cap R[X_1, \ldots, X_n]}$ is a DVR with
maximal ideal $M=(Q \cap R[X_1, \ldots, X_n])V$, and, since
$\op{ht} Q=1$, we must have $Q$ extended from $Q \cap R[X_1,
\ldots, X_n]$.  It is then easy to see that $\rx_Q =
V[\{X_{\al}\}]_{M[\{X_{\al}\}]}$ is a DVR.
\end{proof}

\begin{lemma} \label{l:uppers2} Let $Q$ be an upper to zero in
$\rx$ which is also a maximal $t$-ideal.  Then $\op{ht} Q=1$.
\end{lemma}

\begin{proof} First suppose $\{X_{\al}\}=\{X_1, \ldots, X_n\}$.
The result clearly holds if $n=1$ (even if $Q$ is not a maximal
$t$-ideal!). Suppose $n>1$, and let $q=Q \cap R[X_1, \ldots,
X_{n-1}]$. If $q=0$, then $\op{ht} Q=1$ by the case $n=1$. If $q
\ne 0$, then by \cite[Theorem 1.4]{hz}, $q$ is a maximal $t$-ideal
of $R[X_1, \ldots, X_{n-1}]$, and $Q=q[X_n]$.  By induction
$\op{ht} q=1$, and $V=R[X_1, \ldots, X_{n-1}]_q$ is a DVR by
Lemma~\ref{l:uppers}.  Hence $R[X_1, \ldots, X_n]_Q=V[X_n]_Q$ is
also a DVR, and $\op{ht} Q=1$.

For the general case, we may pick $X_1, \ldots, X_n \in
\{X_{\al}\}$ with $q=Q \cap R[X_1, \ldots, X_n] \ne 0$.  By
\cite[Proposition 2.2]{fgh}, $q$ is a maximal $t$-ideal of $R[X_1,
\ldots, X_n]$, and $Q$ is extended from $q$.  The argument now
proceeds as in the induction step above.
\end{proof}

The following extends \cite[Theorem 1.4 and Corollary 1.5]{hz} to
the case of infinitely many indeterminates.

\begin{theorem} \label{t:uppersmaxt} Let $Q$ be an upper to zero
in $\rx$.  Then following statments are equivalent.
\begin{itemize}
\item[(1)] $Q$ is a maximal $t$-ideal.
\item[(2)] $Q$ is $t$-invertible.
\item[(3)] $c(Q)_t=R$, and $\op{ht} Q = 1$. (In this case, a
standard argument shows that $Q$ contains an element $g$ with
$c(g)_v=R$.)
\end{itemize}
In case these equivalent statements hold, then $Q=fK[\{X_{\al}\}]
\cap \rx$ for some $f \in \rx$ such that $f$ is irreducible in
$K[\{X_{\al}\}]$; moreover, we have $Q=(f,g)_v$.
\end{theorem}

\begin{proof} (3) $\ra$ (1): Since $\op{ht} Q=1$, $Q$ is a
$t$-ideal.  Hence $Q$ is contained in a maximal $t$-ideal, say
$N$.  Since $c(Q)_t=R$, we cannot have $N$ extended from $N \cap
R$, whence $N$ is an upper to zero by \cite[Proposition 2.2]{fgh}.
By Lemma~\ref{l:uppers2}, $\op{ht}N=1$, whence $Q=N$, and $Q$ is a
maximal $t$-ideal.

The proofs of (1) $\ra$ (3) and (2) $\ra$ (1) are as in
\cite[Theorem 1.4]{hz}.

(1) $\ra$ (2): This also goes through essentially as in the proof
of \cite[Theorem 1.4]{hz}.  That proof contains a reference to
\cite[Proposition 1.8]{hhj}, which is stated for the case of one
indeterminate.  However, the proof of this latter result extends
to the case of an arbitrary set of indeterminates.  (The content
formula \cite[Corollary 28.3]{g1} is needed.)

To prove the last statement, note that $Q=fK[\{X_{\al}\}] \cap
\rx$ by Lemmas~\ref{l:uppers} and \ref{l:uppers2}.  The fact
$Q=(f,g)_v$ may be proved as in \cite[Corollary 1.5]{hz}.
\end{proof}

\begin{theorem} \label{t:sharppoly} If $R$ is a $t\#$-domain, then
so is $\rx$.
\end{theorem}

\begin{proof} We wish to show that condition (5) of
Theorem~\ref{t:sharpchar} is satisfied. Thus let $N$ be a maximal
$t$-ideal of $\rx$. By \cite[Proposition 2.2]{fgh}, either $N \cap
R=0$ or $N=(N \cap R)\rx$.  In the former case $N$ is divisorial
(being a $t$-invertible $t$-ideal), and $N$ is certainly not
contained in any other maximal $t$-ideal of $\rx$. In the latter
case, $N \cap R$ contains a divisorial ideal $I$ which is
contained in no other maximal $t$-ideal of $R$, and it follows
that $I\rx$ is a divisorial ideal of $\rx$ which is contained in
$N$ and no other maximal $t$-ideal of $\rx$.
\end{proof}

We have been unable to prove the converse of
Theorem~\ref{t:sharppoly}. (Indeed, we doubt that the converse is
true.) However, we can prove that several standard localizations
of $\rx$ are simultaneously $t\#$. We denote by $\rxx$ the ring of
fractions of $\rx$ with respect to the multiplicatively closed
subset of $\rx$ consisting of the polynomials of unit content.
Finally, if $S=\{f \in \rx \mid c(f)_v=R\}$, we denote by $\rxxx$
the ring $\rx_S$. We then have the following description of the
maximal $t$-ideals in these rings.

\begin{lemma} \label{l:maxt}  Denote by $\mc U_1$ the set of
uppers to zero which are also maximal $t$-ideals in $\rx$ and by
$\mc U_2$ the set of those elements $P \in \mc U_1$ which satisfy
$c(P) \ne R$. Then:
\begin{itemize}
\item[(1)] $t\Max(\rx)=\{M\rx \mid M \in t\Max (R)\}
\bigcup \mc U_1$;
\item[(2)] $t\Max (\rxx)=\{M\rxx \mid M \in t\Max (R)\} \bigcup
\{P\rxx \mid P \in \mc U_2\}$;
\item[(3)] $t\Max (\rxxx)=\{M\rxxx \mid M \in t\Max (R)\}$.
\end{itemize}
\end{lemma}

\begin{proof}  (1) Each maximal $t$-ideal of $\rx$ must have the
form indicated by \cite[Proposition 2.2]{fgh}. The reverse
inclusion follows from \cite[Lemma 2.1]{fgh}.

(2) By \cite[Corollary 2.3]{k}, $M\rxx$ is a $t$-ideal of $\rxx$
for each $M \in \tmax(R)$.  Suppose that for some $N \in
\tmax(\rxx)$ we have $N \supseteq M\rxx$. Then since $\rxx$ is a
ring of fractions of $\rx$, $N$ is extended from a maximal
$t$-ideal of $\rx$, which in turn must be extended from a maximal
$t$-ideal of $R$.  It follows that $N=M\rxx$. Hence $M\rxx \in
\tmax{\rxx}$.  If $P \in \mc U_2$, then, since $c(P) \ne R$,
$P\rxx \ne \rxx$. Moreover, since $\op{ht} P=1$ by
Theorem~\ref{t:uppersmaxt}, $\op{ht} P\rxx = 1$ also, and $P\rxx$
is a $t$-prime of $\rxx$. Any maximal $t$-ideal of $\rxx$
containing $P\rxx$ must be extended from a $t$-prime of $\rx$
containing $P$. Therefore, since $P \in \tmax(\rx)$, $P\rxx \in
\tmax(\rxx)$. That each maximal $t$-ideal of $\rxx$ must be of the
form indicated follows from (1) (and the fact that $\rxx$ is a
ring of fractions of $\rx$).

(3)  This follows from the facts that $\rxxx$ is a localization of
$\rxx$ and that each $P \in U_1$ satisfies $c(P)_t = \rx$ by
Theorem~\ref{t:uppersmaxt} so that $P\rxxx = \rxxx$.
\end{proof}

In the proof of the following result, we often invoke
Lemma~\ref{l:maxt} without explicit reference.

\begin{theorem} \label{t:sharppoly2} The following statements are
equivalent for a domain $R$. \begin{itemize}

\item[(1)] $\rx$ is a $t\#$-domain.

\item[(2)] $R(\{X_{\al}\})$ is a $t\#$-domain.

\item[(3)] $\rxxx$ is a $t\#$-domain.

\end{itemize}

If, in addition, $\rx$ is $v$-coherent, then these conditions are
equivalent to: $R$ is a $t\#$-domain.
\end{theorem}

\begin{proof} (1) $\ra$ (2): If $M \in \tmax(R)$, then by
Theorem~\ref{t:sharpchar}, $M$ contains a divisorial ideal $I$
which is contained in no other maximal $t$-ideal of $R$.  It
follows that $I\rxx$ is a divisorial ideal of $\rxx$
\cite[Corollary 2.3]{k}, and it is clear that $I\rxx$ is contained
in $M\rxx$ but in no other maximal $t$-ideal of $\rxx$. Hence each
maximal $t$-ideal of $\rxx$ of the form $M\rxx$ contains a
divisorial ideal contained in no other maximal $t$-ideal of
$\rxx$.  On the other hand, if $P\rxx$ is a maximal $t$-ideal of
$\rxx$ with $P \in \mc U_2$, then $P$ is divisorial, from which it
follows $P\rxx$ is also divisorial (and is clearly not contained
in any other maximal $t$-ideal of $\rxx$). By
Theorem~\ref{t:sharpchar}, $\rxx$ is a $t\#$-domain.

(2) $\ra$ (3): Similar (but easier).

(3) $\ra$ (1): Let $M$ be a maximal $t$-ideal of $R$.  By
hypothesis and Theorem~\ref{t:sharpchar} ((1) $\lra$ (6)), there
is an element $u \in K(\{X_{\alpha}\})$ such that $(\rxxx:_{\rxxx}
u)$ is contained in $M\rxxx$ and no other maximal $t$-ideal of
$\rxxx$.  Let $I=(\rx:_{\rx} u)$.  Then $I$ is divisorial in
$\rx$, and $I\rxxx = (\rxxx:_{\rxxx} u)$.  Clearly, $I \sub M\rx$
and $I \nsubseteq N\rx$ for each maximal $t$-ideal $N$ of $R$ with
$N \ne M$.  Moreover, $I$ is contained in at most finitely many
maximal $t$-ideals $P$ with $P \cap R=0$.  We shall show how to
enlarge $I$ so as to avoid each such $P$.  By
Theorem~\ref{t:uppersmaxt}, we have that $P\rxxx=\rxxx$, and $P$
is $v$-finite. Therefore, since $\rx_P$ is a DVR, we may pick $h
\in \rx \sm P$ with $hP^n \sub I$.  Hence $h\rxxx=hP^n\rxxx \sub
I\rxxx$, and there is an element $g \in \rx$ with $c(g)_v=R$ and
$gh \in I$.  In particular, $g \notin M\rx$, so that the
divisorial ideal $(I:_{\rx} g)$ is contained in $M\rx$. Moreover,
$h \in (I:_{\rx} g) \sm P$.  Hence $(I:_{\rx} g)$ is a divisorial
ideal contained in $M\rx \sm P$.  This process may be continued
finitely many times to produce a divisorial ideal which is
contained in $M\rx$ but in no other maximal $t$-ideal of $\rx$.
Thus $\rx$ is a $t\#$-domain.

To prove the final statement, assume that $\rx$ is a $t\#$-domain.
Let $M \in \tmax (R)$.  By Theorem~\ref{t:sharpchar2}, there is a
finitely generated ideal $I$ of $\rx$ such that $I \sub M\rx$ and
$I$ is contained in no other maximal $t$-ideal of $\rx$.  Clearly,
$c(I) \sub M$ and $c(I)$ is contained in no other maximal
$t$-ideal of $R$.  Another application of
Theorem~\ref{t:sharpchar2} completes the proof.
\end{proof}

It is well known that a domain $R$ is a PVMD if and only if $\rx$
is a PVMD.  Thus for a PVMD $R$ the conditions of
Theorem~\ref{t:sharppoly2} are each equivalent to $R$ being a
$t\#$-domain.  It follows that if $R$ is a Pr\"ufer domain, then
$R$ is a $\#$-domain if and only if $\rx$ is a $t\#$-domain.

Now recall that it is possible for a polynomial ring over a Mori
domain to fail to be Mori \cite{r}.  In view of the fact that a
Mori domain is automatically a $t\#$-domain
(Corollary~\ref{c:tv}), we see by Theorem~\ref{t:sharppoly} that
if $R$ is a Mori domain, then $\rx$ is a $t\#$-domain even though
$\rx$ may not be a Mori domain.

It is an open question whether $R$ $v$-coherent implies that $\rx$
is $v$-coherent.  We are therefore unable to determine whether the
last statement of Theorem~\ref{t:sharppoly2} remains true if we
assume only that $R$ is $v$-coherent.  It is true, however, that
$v$-coherence of $\rx$ implies that of $R$, as the following
result shows.

\begin{proposition} \label{p:vcohpoly} If $\rx$ is $v$-coherent,
then $R$ is $v$-coherent.
\end{proposition}

\begin{proof} Let $I$ be a finitely generated ideal of $R$.  We
have $(I[\{X_{\al}\}])^{-1}=I^{-1}[\{X_{\al}\}]$; by hypothesis,
this produces a finitely generated fractional ideal $J$ of $\rx$
with $I^{-1}[\{X_{\al}\}]=J_v$.  We may assume $1 \in J$.
Moreover, since $\rx \sub I^{-1}[\{X_{\al}\}] \sub
K[\{X_{\al}\}]$, we have $\rx \sub J \sub K[\{X_{\al}\}]$.  Hence
$c(J)$ is a finitely generated ideal of $R$ with $1 \in c(J)$.  We
claim that $c(J)_v=I^{-1}$.  Note that $J \sub c(J)[\{X_{\al}\}]
\sub I^{-1}[\{X_{\al}\}]$.  Hence $$J_v \sub
(c(J))[\{X_{\al}\}])_v = c(J)_v[\{X_{\al}\}] \sub
I^{-1}\{X_{\al}\}] = J_v,$$ and the claim follows.
\end{proof}




\end{document}